\documentclass[12pt]{article}
\usepackage{amsmath}
\usepackage{amsfonts}
\usepackage{amssymb}

\def\thebibliograph#1#2{\section*{{\normalsize \bf #2}}\list
   {[\arabic{enumi}]}{\settowidth\labelwidth{[#1]}\leftmargin\labelwidth
     \advance\leftmargin\labelsep
     \usecounter{enumi}}
     \def\newblock{\hskip .11em plus .33em minus -.07em}
     \sloppy
     \sfcode`\.=1000\relax}

\newtheorem{theorem}{Theorem}

\newtheorem{definition}{Definition}
\newtheorem{corollary}{Corollary}
\newtheorem{lemma}{Lemma}
\newtheorem{remark}{Remark}
\input{tcilatex}
\begin{document}

\title{Sobolev embeddings for Herz-type Triebel-Lizorkin spaces}
\author{ Douadi Drihem \ \thanks{%
M'sila University, Department of Mathematics, Laboratory of Functional
Analysis and Geometry of Spaces , P.O. Box 166, M'sila 28000, Algeria,
e-mail: \texttt{\ douadidr@yahoo.fr}}}
\date{\today }
\maketitle

\begin{abstract}
In this paper we prove the Sobolev embeddings for Herz-type Triebel-Lizorkin
spaces, 
\begin{equation*}
\dot{K}_{q}^{\alpha _{2},r}F_{\theta }^{s_{2}}\hookrightarrow \dot{K}%
_{s}^{\alpha _{1},p}F_{\beta }^{s_{1}}
\end{equation*}%
where the parameters $\alpha _{1},\alpha _{2},s_{1},s_{2},s,q,r,p,\beta $
and $\theta $ satisfy some suitable conditions. An application we obtain new
embeddings between Herz and Triebel-Lizorkin spaces. Moreover, we present
the Sobolev embeddings for Triebel-Lizorkin spaces equipped with power
weights. All these results cover the results on classical Triebel-Lizorkin
spaces. \medskip


\end{abstract}

\section{Introduction}


Function spaces have been widely used in various areas of analysis such as
harmonic analysis and partial differential equations. In recent years, there
has been increasing interest in a new family of function spaces which
generalize the Besov spaces and Triebel-Lizorkin spaces. Some example of
these spaces can be mentioned such as Herz-type Triebel-Lizorkin spaces, $%
\dot{K}_{q}^{\alpha ,p}F_{\beta }^{s}$, that initially appeared in the
papers of J. Xu and D. Yang \cite{XuYang03} , \cite{XuYang05} and \cite{Xu03}%
. Several basic properties were established, such as the Fourier analytical
characterisation and\ lifting properties. When $\alpha =0$ and$\ p=q$ they
coincide with the usual function spaces $F_{p,q}^{s}$.\vskip5pt

The interest in these spaces comes not only from theoretical reasons but
also from their applications to several classical problems in analysis. In 
\cite{LuYang97}, Lu and Yang introduced the Herz-type Sobolev and Bessel
potential spaces. They gave some applications to partial differential
equations. Also in \cite{T11}, Y. Tsutsui, studied the Cauchy problem for
Navier-Stokes equations on Herz spaces\ and weak Herz spaces.\vskip5pt

Since the Sobolev embedding plays an important role in theory of function
spaces and PDE's, the main aim of this paper is to prove the Sobolev
embedding of $\dot{K}_{q}^{\alpha ,p}F_{\beta }^{s}$ spaces. First we shall
prove the Sobolev embeddings of associated sequence spaces. Then, from the
so-called $\varphi $-transform characterization in the sense of Frazier and
Jawerth, we deduce the main result of this paper. As a consequence, we
obtain new Jawerth-Franke-type embeddings, the Sobolev embeddings for
Triebel-Lizorkin spaces equipped with power weights, new embeddings between
Herz and Triebel-Lizorkin spaces, and we present some remarks about the
wavelet caracterization of Herz-Triebel-Lizorkin spaces. All these results
generalize the existing classical results on Triebel-Lizorkin spaces.

To recall the definition of these function spaces, we need some notation.
For any $u>0,k\in \mathbb{Z}$ we set $C\left( u\right) =\{x\in \mathbb{R}%
^{n}:u/2\leq \left\vert x\right\vert <u\}$ and $C_{k}=C(2^{k})$. For $x\in 
\mathbb{R}^{n}$ and $r>0$ we denote by $B(x,r)$ the open ball in $\mathbb{R}%
^{n}$ with center $x$ and radius $r$. Let $\chi _{k}$, for $k\in \mathbb{Z}$%
, denote the characteristic function of the set $C_{k}$. The expression $%
f\approx g$ means that $C$ $g\leq f\leq c\,g$ for some independent constants 
$c,C$ and non-negative functions $f$ and $g$. \vskip5pt

\noindent We denote by $\left\vert \Omega \right\vert $ the $n$-dimensional
Lebesgue measure of $\Omega \subseteq \mathbb{R}^{n}$. For any measurable
subset $\Omega \subseteq \mathbb{R}^{n}$ the Lebesgue space $L^{p}(\Omega )$%
, $0<p\leq \infty $ consists of all measurable functions for which $%
\left\Vert f\mid L^{p}(\Omega )\right\Vert =\left( \int_{\Omega }\left\vert
f(x)\right\vert ^{p}dx\right) ^{1/p}<\infty $, $0<p<\infty $ and $\left\Vert
f\mid L^{\infty }(\Omega )\right\Vert =\underset{x\in \Omega }{\text{ess-sup}%
}\left\vert f(x)\right\vert <\infty $. If $\Omega =\mathbb{R}^{n}$ we put $%
L^{p}(\mathbb{R}^{n})=L^{p}$ and $\left\Vert f\mid L^{p}(\mathbb{R}%
^{n})\right\Vert =\left\Vert f\right\Vert _{p}$. The Hardy-Littlewood
maximal operator $\mathcal{M}$ is defined on locally integrable functions by 
\begin{equation*}
\mathcal{M}f(x)=\sup_{r>0}\frac{1}{\left\vert B(x,r)\right\vert }%
\int_{B(x,r)}\left\vert f\left( y\right) \right\vert dy
\end{equation*}%
and $\mathcal{M}_{t}f=\left( \mathcal{M}|f|^{t}\right) ^{1/t}$ for any $%
0<t\leq 1$. The symbol $\mathcal{S}(\mathbb{R}^{n})$ is used in place of the
set of all Schwartz functions $\varphi $ on $\mathbb{R}^{n}$ and we denote
by $\mathcal{S}^{\prime }(\mathbb{R}^{n})$ the dual space of all tempered
distributions on $\mathbb{R}^{n}$. We define the Fourier transform of a
function $f\in \mathcal{S}(\mathbb{R}^{n})$ by $\mathcal{F(}f)(\xi )=\left(
2\pi \right) ^{-n/2}\int_{\mathbb{R}^{n}}e^{-ix\cdot \xi }f(x)dx$. Its
inverse is denoted by $\mathcal{F}^{-1}f$. Both $\mathcal{F}$ and $\mathcal{F%
}^{-1}$ are extended to the dual Schwartz space $\mathcal{S}^{\prime }(%
\mathbb{R}^{n})$ in the usual way.\vskip5pt

Let $\mathbb{Z}^{n}$ be the lattice of all points in $\mathbb{R}^{n}$ with
integer-valued components. If $v\in \mathbb{N}_{0}$ and $m=(m_{1},...,m_{n})%
\in \mathbb{Z}^{n}$ we denote $Q_{v,m}$ the dyadic cube in $\mathbb{R}^{n}$%
\begin{equation*}
Q_{v,m}=\{(x_{1},...,x_{n}):m_{i}\leq 2^{v}x_{i}<m_{i}+1,i=1,2,...,n\}.
\end{equation*}%
By $\chi _{v,m}$ we denote the characteristic function of the cube $Q_{v,m}$.

Given two quasi-Banach spaces $X$ and $Y$, we write $X\hookrightarrow Y$ if $%
X\subset Y$ and the natural embedding of $X$ in $Y$ is continuous. We use $c$
as a generic positive constant, i.e.\ a constant whose value may change from
appearance to appearance.\vskip5pt

\section{Function spaces}

We start by recalling the definition and some of the properties of the
homogenous Herz spaces $\dot{K}_{q}^{\alpha ,p}$.\vskip5pt

\begin{definition}
\label{Herz-spaces} \textit{Let }$\alpha \in \mathbb{R},0<p,q\leq \infty $%
\textit{. The homogeneous Herz space }$\dot{K}_{q}^{\alpha ,p}$\textit{\ is
defined by }%
\begin{equation*}
\dot{K}_{q}^{\alpha ,p}=\{f\in L_{\mathrm{loc}}^{q}(\mathbb{R}^{n}\setminus
\{0\}):\big\|f\mid \dot{K}_{q}^{\alpha ,p}\big\|<\infty \},
\end{equation*}%
\textit{where }%
\begin{equation*}
\big\|f\mid \dot{K}_{q}^{\alpha ,p}\big\|=\Big(\sum\limits_{k=-\infty
}^{\infty }2^{k\alpha p}\text{ }\left\Vert f\chi _{k}\right\Vert _{q}^{p}%
\Big)^{1/p},
\end{equation*}%
\textit{with the usual modifications made when }$p=\infty $\textit{\ and/or }%
$q=\infty $\textit{.}
\end{definition}

The spaces $\dot{K}_{q}^{\alpha ,p}$ are quasi-Banach spaces and if $\min
(p,q)\geq 1$ then $\dot{K}_{q}^{\alpha ,p}$ are Banach spaces. When $\alpha
=0$ and $0<p=q\leq \infty $ then $\dot{K}_{p}^{0,p}$ coincides with the
Lebesgue spaces $L^{p}$. Various important results have been proved in the
space $\dot{K}_{q}^{\alpha ,p}$ under some assumptions on $\alpha ,p$ and $q$%
. The conditions $-\frac{n}{q}<\alpha <n(1-\frac{1}{q}),1<q<\infty $ and $%
0<p\leq \infty $ is crucial in the study of the boundedness of classical
operators in $\dot{K}_{q}^{\alpha ,p}$ spaces. This fact was first realized
by Li and Yang \cite{LiYang96} with the proof of the boundedness of the
maximal function. The proof of the main result of this section is based on
the following result, see Tang and Yang \cite{TD00}.

\begin{lemma}
\label{Maximal-Inq}Let $1<\beta <\infty ,1<q<\infty $ and $0<p\leq \infty $.
If $\{f_{j}\}_{j=0}^{\infty }$ is a sequence of locally integrable functions
on $\mathbb{R}^{n}$ and $-\frac{n}{q}<\alpha <n(1-\frac{1}{q})$, then%
\begin{equation*}
\Big\|\Big(\sum_{j=0}^{\infty }(\mathcal{M}f_{j})^{\beta }\Big)^{1/\beta }|%
\dot{K}_{q}^{\alpha ,p}\Big\|\leq c\Big\|\Big(\sum_{j=0}^{\infty
}|f_{j}|^{\beta }\Big)^{1/\beta }|\dot{K}_{q}^{\alpha ,p}\Big\|.
\end{equation*}
\end{lemma}

A detailed discussion of the properties of these spaces my be found in the
papers \cite{HerYang99}, \cite{LuYang1.95}, \cite{LuYang2.95}, and
references therein.

\ \ \ Now, we\ present the Fourier analytical definition of Herz-type
Triebel-Lizorkin spaces\textit{\ }$\dot{K}_{q}^{\alpha ,p}F_{\beta }^{s}$
and recall their basic properties. We first need the concept of a smooth
dyadic resolution of unity. Let $\phi _{0}$\ be a function\ in $\mathcal{S}(%
\mathbb{R}^{n})$\ satisfying $\phi _{0}(x)=1$\ for\ $\left\vert x\right\vert
\leq 1$\ and\ $\phi _{0}(x)=0$\ for\ $\left\vert x\right\vert \geq 2$.\ We
put $\phi _{j}(x)=\phi _{0}(2^{-j}x)-\phi _{0}(2^{1-j}x)$ for $j=1,2,3,....$
Then $\{\phi _{j}\}_{j\in \mathbb{N}_{0}}$\ is a resolution of unity, $%
\sum_{j=0}^{\infty }\phi _{j}(x)=1$ for all $x\in \mathbb{R}^{n}$.\ Thus we
obtain the Littlewood-Paley decomposition $f=\sum_{j=0}^{\infty }\mathcal{F}%
^{-1}\phi _{j}\ast f$ of all $f\in \mathcal{S}^{\prime }(\mathbb{R}^{n})$ $($%
convergence in $\mathcal{S}^{\prime }(\mathbb{R}^{n}))$.

We are now in a position to state the definition of Herz-type
Triebel-Lizorkin spaces.

\begin{definition}
\label{Herz-Besov-Triebel}\textit{Let }$\alpha ,s\in \mathbb{R},0<p,q<\infty 
$\textit{\ and }$0<\beta \leq \infty $\textit{. The Herz-type
Triebel-Lizorkin space }$\dot{K}_{q}^{\alpha ,p}F_{\beta }^{s}$ \textit{is
the collection of all} $f\in \mathcal{S}^{\prime }(\mathbb{R}^{n})$\textit{\
such that}%
\begin{equation}
\big\|f\mid \dot{K}_{q}^{\alpha ,p}F_{\beta }^{s}\big\|=\Big\|\Big(%
\sum\limits_{j=0}^{\infty }2^{js\beta }\left\vert \mathcal{F}^{-1}\phi
_{j}\ast f\right\vert ^{\beta }\Big)^{1/\beta }\mid \dot{K}_{q}^{\alpha ,p}%
\Big\|<\infty ,  \label{H-T-L-space}
\end{equation}%
\textit{with the obvious modification if }$\beta =\infty .$
\end{definition}

\begin{remark}
$\mathrm{Let}$\textrm{\ }$s\in \mathbb{R},0<p,q<\infty ,0<\beta \leq \infty $
$\mathrm{and}$ $\alpha >-n/q$\textrm{. }$\mathrm{The}$ $\mathrm{spaces}$%
\textrm{\ }$\dot{K}_{q}^{\alpha ,p}F_{\beta }^{s}$\textrm{\ }$\mathrm{are}$ $%
\mathrm{independent}$ $\mathrm{of}$ $\mathrm{the}$ $\mathrm{particular}$ $%
\mathrm{choice}$ $\mathrm{of}$ $\mathrm{the}$ $\mathrm{smooth}$ $\mathrm{%
dyadic}$ $\mathrm{resolution}$ $\mathrm{of}$ $\mathrm{unity}$\textrm{\ }$%
\{\phi _{j}\}_{j\in \mathbb{N}_{0}}$ $\mathrm{(in}$ $\mathrm{the}$ $\mathrm{%
sense}$ $\mathrm{of\ equivalent}$ $\mathrm{quasi}$\textrm{-}$\mathrm{norms)}$%
\textrm{. }$\mathrm{In}$ $\mathrm{particular}$ $\dot{K}_{q}^{\alpha
,p}F_{\beta }^{s}$\textrm{\ }$\mathrm{are}$ $\mathrm{quasi}$\textrm{-}$%
\mathrm{Banach}$ $\mathrm{spaces}$ $\mathrm{and}$ $\mathrm{if}$\textrm{\ }$%
p,q,\beta \geq 1\mathrm{,}$ $\mathrm{then}$ $\dot{K}_{q}^{\alpha ,p}F_{\beta
}^{s}$ $\mathrm{are}$ $\mathrm{Banach}$ $\mathrm{spaces}$\textrm{. }$\mathrm{%
Further}$ $\mathrm{results,concerning,}$ $\mathrm{for}$ $\mathrm{instance,}$ 
$\mathrm{lifting}$ $\mathrm{properties,}$ $\mathrm{Fourier}$ $\mathrm{%
multiplier}$ $\mathrm{and}$ $\mathrm{local}$ $\mathrm{means}$ $\mathrm{%
characterizations}$ $\mathrm{can}$ $\mathrm{be}$ $\mathrm{found}$ $\mathrm{in%
}$ $\mathrm{\cite{XuYang05},}$ $\mathrm{\cite{Xu03},}$ $\mathrm{\cite{Xu04}}$
$\mathrm{and}$ $\mathrm{\cite{Xu09}.}$
\end{remark}

Now we give the definitions of the spaces $B_{p,\beta }^{s}$ and $F_{p,\beta
}^{s}$.\vskip5pt

\begin{definition}
$\mathrm{(i)}$\textit{\ Let }$s\in \mathbb{R}$\textit{\ and }$0<p,\beta \leq
\infty $\textit{. The Besov space }$B_{p,\beta }^{s}$ \textit{is the
collection of all} $f\in \mathcal{S}^{\prime }(\mathbb{R}^{n})$\textit{\
such that} 
\begin{equation*}
\left\Vert f\mid B_{p,\beta }^{s}\right\Vert =\Big(\sum\limits_{j=0}^{\infty
}2^{js\beta }\left\Vert \mathcal{F}^{-1}\phi _{j}\ast f\right\Vert
_{p}^{\beta }\Big)^{1/\beta }<\infty .
\end{equation*}%
$\mathrm{(ii)}$\textit{\ Let }$s\in \mathbb{R},0<p<\infty $\textit{\ and }$%
0<\beta \leq \infty $\textit{. The Triebel-Lizorkin space }$F_{p,\beta }^{s}$
\textit{is the collection of all} $f\in \mathcal{S}^{\prime }(\mathbb{R}%
^{n}) $\textit{\ such that}%
\begin{equation*}
\left\Vert f\mid F_{p,\beta }^{s}\right\Vert =\Big\|\Big(\sum\limits_{j=0}^{%
\infty }2^{js\beta }\left\vert \mathcal{F}^{-1}\phi _{j}\ast f\right\vert
^{\beta }\Big)^{1/\beta }\Big\|_{p}<\infty .
\end{equation*}
\end{definition}

The theory of the spaces $B_{p,\beta }^{s}$ and $F_{p,\beta }^{s}$ has been
developed in detail in \cite{Triebel83}, \cite{Triebel92} and \cite%
{Triebel06} but has a longer history already including many contributors; we
do not want to discuss this here. Clearly, for $s\in \mathbb{R},0<p<\infty $
and $0<\beta \leq \infty ,$%
\begin{equation*}
\dot{K}_{p}^{0,p}F_{\beta }^{s}=F_{p,\beta }^{s}.
\end{equation*}

Let us consider $k_{0},k\in \mathcal{S}(\mathbb{R}^{n})$ and $S\geq -1$ an
integer such that for an $\varepsilon >0$%
\begin{align}
\left\vert \mathcal{F}k_{0}(\xi )\right\vert & >0\text{\quad for\quad }%
\left\vert \xi \right\vert <2\varepsilon  \label{T-cond1} \\
\left\vert \mathcal{F}k(\xi )\right\vert & >0\text{\quad for\quad }\frac{%
\varepsilon }{2}<\left\vert \xi \right\vert <2\varepsilon  \label{T-cond2}
\end{align}%
and%
\begin{equation}
\int_{\mathbb{R}^{n}}x^{\alpha }k(x)dx=0\text{\quad for any\quad }\left\vert
\alpha \right\vert \leq S.  \label{moment-cond}
\end{equation}%
Here $\mathrm{\eqref{T-cond1}}$ and $\mathrm{\eqref{T-cond2}}$\ are
Tauberian conditions, while $\mathrm{\eqref{moment-cond}}$ are moment
conditions on $k$. We recall the notation%
\begin{equation*}
k_{t}(x)=t^{-n}k(t^{-1}x)\text{,}\quad k_{j}(x)=k_{2^{-j}}(x)\text{,\quad for%
}\quad t>0\quad \text{and}\quad j\in \mathbb{N}.
\end{equation*}%
Usually $k_{j}\ast f$ is called local mean. The following result is from 
\cite[Theorem 1]{Xu05}.

\begin{theorem}
\label{loc-mean-char}\textit{Let }$\alpha ,s\in \mathbb{R},0<p,q<\infty $, $%
0<\beta \leq \infty $, $\alpha >-n/q$ and $s<S+1$\textit{. Then}%
\begin{equation}
\big\|f\mid \dot{K}_{q}^{\alpha ,p}F_{\beta }^{s}\big\|^{\prime }=\Big\|%
\left( 2^{js}k_{j}\ast f\right) _{j}\mid \dot{K}_{q}^{\alpha ,p}(\ell
_{\beta })\Big\|,  \label{equiv-norm3}
\end{equation}%
\textit{is equivalent quasi-norm on }$\dot{K}_{q}^{\alpha ,p}F_{\beta }^{s}$.
\end{theorem}

We introduce the sequence spaces associated with the function spaces $\dot{K}%
_{q}^{\alpha ,p}F_{\beta }^{s}$. If 
\begin{equation*}
\lambda =\{\lambda _{v,m}\in \mathbb{C}:v\in \mathbb{N}_{0},m\in \mathbb{Z}%
^{n}\},
\end{equation*}%
$\alpha ,s\in \mathbb{R},0<p,q<\infty $\ and $0<\beta \leq \infty $, we set%
\begin{equation}
\big\|\lambda \mid \dot{K}_{q}^{\alpha ,p}f_{\beta }^{s}\big\|=\Big\|\Big(%
\sum_{v=0}^{\infty }\sum\limits_{m\in \mathbb{Z}^{n}}2^{vs\beta }|\lambda
_{v,m}|^{\beta }\chi _{v,m}\Big)^{1/\beta }\mid \dot{K}_{q}^{\alpha ,p}\Big\|%
.  \label{atomic-norm}
\end{equation}%
Let $\Phi ,\psi ,\varphi $ and $\Psi $ satisfy 
\begin{equation}
\Phi ,\Psi ,\varphi ,\psi \in \mathcal{S}(\mathbb{R}^{n})  \label{Ass1}
\end{equation}%
\begin{equation}
\text{supp}\mathcal{F}\Phi \text{, supp}\mathcal{F}\Psi \subset \overline{%
B(0,2)}\text{ such that }|\mathcal{F}\Phi (\xi )|,|\mathcal{F}\Psi (\xi
)|\geq c\text{ if }|\xi |\leq \frac{5}{3}  \label{Ass2}
\end{equation}%
and 
\begin{equation}
\text{supp}\mathcal{F}\varphi \text{, supp}\mathcal{F}\psi \subset \overline{%
B(0,2)}\backslash B(0,1/2)\text{ such that }|\mathcal{F}\varphi (\xi )|,|%
\mathcal{F}\psi (\xi )|\geq c\text{ if }\frac{3}{5}\leq |\xi |\leq \frac{5}{3%
}  \label{Ass3}
\end{equation}%
such that%
\begin{equation}
\overline{\mathcal{F}\Phi (-\xi )}\mathcal{F}\Psi (\xi )+\sum_{j=1}^{\infty }%
\overline{\mathcal{F}\varphi (-2^{-j}\xi )}\mathcal{F}\psi (2^{-j}\xi
)=1,\quad \xi \in \mathbb{R}^{n},  \label{Ass4}
\end{equation}%
where $c>0$. Recall that the $\varphi $-transform $S_{\varphi }$ is defined
by setting $(S_{\varphi }f)_{0,m}=\langle f,\Phi _{m}\rangle $ where $\Phi
_{m}(x)=\Phi (x-m)$ and $(S_{\varphi }f)_{v,m}=\langle f,\varphi
_{v,m}\rangle $ where $\varphi _{v,m}(x)=2^{vn/2}\varphi (2^{v}x-m)$ and $%
v\in \mathbb{N}$. The inverse $\varphi $-transform $T_{\psi }$ is defined by 
\begin{equation*}
T_{\psi }\lambda =\sum_{m\in \mathbb{Z}^{n}}\lambda _{0,m}\Psi
_{m}+\sum_{v=1}^{\infty }\sum_{m\in \mathbb{Z}^{n}}\lambda _{v,m}\psi _{v,m},
\end{equation*}%
where $\lambda =\{\lambda _{v,m}\in \mathbb{C}:v\in \mathbb{N}_{0},m\in 
\mathbb{Z}^{n}\}$, see \cite{FJ90}.

For a sequence $\lambda =\{\lambda _{v,m}\in \mathbb{C}:v\in \mathbb{N}%
_{0},m\in \mathbb{Z}^{n}\},0<r\leq \infty $ and a fixed $d>0$, set%
\begin{equation*}
\lambda _{v,m,r,d}^{\ast }=\Big(\sum_{h\in \mathbb{Z}^{n}}\frac{|\lambda
_{v,h}|^{r}}{(1+2^{v}|2^{-v}h-2^{-v}m|)^{d}}\Big)^{1/r}
\end{equation*}%
and $\lambda _{r,d}^{\ast }=\{\lambda _{v,m,r,d}^{\ast }\in \mathbb{C}:v\in 
\mathbb{N}_{0},m\in \mathbb{Z}^{n}\}$.

\begin{lemma}
\label{lamda-equi}\textit{Let }$\alpha ,s\in \mathbb{R},0<p\leq \infty
,0<q<\infty ,0<\beta \leq \infty $, $d>n$ \textit{and }$\alpha >-n/q$. Then%
\begin{equation*}
\Big\|\lambda _{\min (q,\frac{n}{\frac{n}{q}+\alpha },\beta ),d}^{\ast }\mid 
\dot{K}_{q}^{\alpha ,p}f_{\beta }^{s}\Big\|\approx \big\|\lambda \mid \dot{K}%
_{q}^{\alpha ,p}f_{\beta }^{s}\big\|
\end{equation*}
\end{lemma}

\textbf{Proof.} Obviously, 
\begin{equation*}
\left\Vert \lambda \mid \dot{K}_{q}^{\alpha ,p}f_{\beta }^{s}\right\Vert
\leq \Big\|\lambda _{\min (q,\frac{n}{\frac{n}{q}+\alpha },\beta ),d}^{\ast
}\mid \dot{K}_{q}^{\alpha ,p}f_{\beta }^{s}\Big\|.
\end{equation*}%
From Lemma A.2 of \cite{FJ90}, we obtain%
\begin{equation*}
\lambda _{v,m,\min (q,\frac{n}{\frac{n}{q}+\alpha },\beta ),d}^{\ast }\leq c%
\mathcal{M}_{a}\Big(\sum_{h\in \mathbb{Z}^{n}}|\lambda _{v,h}|\chi _{v,h}%
\Big)(x),\quad x\in Q_{v,m},
\end{equation*}%
where, 
\begin{equation*}
0<a\leq r=\min (q,\frac{n}{\frac{n}{q}+\alpha },\beta )<\infty ,\text{ }da>nr
\end{equation*}%
and $c>0$ depend only $n$ and $d$. Let $\varepsilon =\frac{d}{n}-1>0$ and $a=%
\tfrac{r}{1+\varepsilon /2}$, then $0<a<r$ and $da>nr$. Hence 
\begin{eqnarray*}
&&\Big\|\lambda _{\min (q,\frac{n}{\frac{n}{q}+\alpha },\beta ),d}^{\ast
}\mid \dot{K}_{q}^{\alpha ,p}f_{\beta }^{s}\Big\| \\
&\leq &c\Big\|\Big(\sum_{v=0}^{\infty }\mathcal{M}_{a/\beta }\Big(\sum_{h\in 
\mathbb{Z}^{n}}2^{vs\beta }|\lambda _{v,h}|^{\beta }\chi _{v,h}\Big)\Big)%
^{a/\beta }\mid \dot{K}_{q/a}^{\alpha a,p/a}\Big\|^{1/a}.
\end{eqnarray*}%
Observe that $\frac{\beta }{a}>1$, $\frac{q}{a}>1$ and $\frac{-na}{q}<\alpha
a<n(1-\frac{a}{q})$. Applying Lemma \ref{Maximal-Inq} to estimate the last
expression by%
\begin{equation*}
c\Big\|\Big(\sum_{v=0}^{\infty }\sum_{h\in \mathbb{Z}^{n}}2^{vs\beta
}|\lambda _{v,h}|^{\beta }\chi _{v,h}\Big)^{1/\beta }\mid \dot{K}%
_{q}^{\alpha ,p}\Big\|=c\big\|\lambda \mid \dot{K}_{q}^{\alpha ,p}f_{\beta
}^{s}\big\|.
\end{equation*}%
The proof of the lemma is thus complete. \ \ \rule{3mm}{3mm}

To prove the main results of this paper we need the following theorem.

\begin{theorem}
\label{phi-tran}\textit{Let }$\alpha ,s\in \mathbb{R},0<p<\infty ,0<q<\infty
,0<\beta \leq \infty $ \textit{and }$\alpha >-n/q$. \textit{Suppose that }$%
\varphi $ and $\Phi $ satisfy \eqref{Ass1}-\eqref{Ass4}. The operators $%
S_{\varphi }:\dot{K}_{q}^{\alpha ,p}F_{\beta }^{s}\rightarrow \dot{K}%
_{q}^{\alpha ,p}f_{\beta }^{s}$ and $T_{\psi }:\dot{K}_{q}^{\alpha
,p}f_{\beta }^{s}\rightarrow \dot{K}_{q}^{\alpha ,p}F_{\beta }^{s}$ are
bounded. Furthermore, $T_{\psi }\circ S_{\varphi }$ is the identity on $\dot{%
K}_{q}^{\alpha ,p}F_{\beta }^{s}$.
\end{theorem}

\textbf{Proof.} We use the same arguments of \cite[Theorem 2.2]{FJ90}, see
also \cite[Theorem 2.1]{SiYY10} and \cite[Theorem 3.12]{Bow07}. For any $%
f\in \mathcal{S}^{\prime }(\mathbb{R}^{n})$ we put $\sup
(f)=\{\sup_{v,m}(f):v\in \mathbb{N}_{0},m\in \mathbb{Z}^{n}\}$ where 
\begin{equation*}
\sup_{v,m}(f)=\sup_{y\in Q_{v,m}}|\widetilde{\varphi _{v}}\ast f(y)|
\end{equation*}%
if $v\in \mathbb{N},m\in \mathbb{Z}^{n}$ and 
\begin{equation*}
\sup_{0,m}(f)=\sup_{y\in Q_{0,m}}|\widetilde{\Phi }\ast f(y)|
\end{equation*}%
if $m\in \mathbb{Z}^{n}$. For any $\gamma \in \mathbb{Z}_{+}$, we define the
sequence $\inf_{\gamma }(f)=\{\inf_{v,m,\gamma }(f):v\in \mathbb{N}_{0},m\in 
\mathbb{Z}^{n}\}$ by setting 
\begin{equation*}
\inf_{v,m,\gamma }(f)=\sup_{h\in \mathbb{Z}^{n}}\{\inf_{y\in Q_{v+\gamma
,h}}|\widetilde{\varphi _{v}}\ast f(y)|:Q_{v+\gamma ,h}\cap Q_{v,m}\neq
\emptyset \}
\end{equation*}%
if $v\in \mathbb{N},m\in \mathbb{Z}^{n}$ and 
\begin{equation*}
\inf_{0,m,\gamma }(f)=\sup_{h\in \mathbb{Z}^{n}}\{\inf_{y\in Q_{\gamma ,h}}|%
\widetilde{\Phi }\ast f(y)|:Q_{\gamma ,h}\cap Q_{0,m}\neq \emptyset \}
\end{equation*}%
if $m\in \mathbb{Z}^{n}$. Here $\widetilde{\varphi _{j}}(x)=2^{jn}\overline{%
\varphi (-2^{j}x)}$ and $\widetilde{\Phi }(x)=\overline{\Phi (-x)}$. As in
Lemma A.5 of \cite{FJ90}, see also \cite{Bow07} and \cite{SiYY10} we obtain 
\begin{equation*}
\big\|\text{inf}_{\gamma }(f)\mid \dot{K}_{q}^{\alpha ,p}f_{\beta }^{s}\big\|%
\leq c\big\|f\mid \dot{K}_{q}^{\alpha ,p}F_{\beta }^{s}\big\|
\end{equation*}%
for any $s\in \mathbb{R},0<p,q<\infty ,0<\beta \leq \infty $, $\alpha >-n/q$
and $\gamma >0$ sufficiently large. Indeed, we have%
\begin{equation*}
\big\|\text{inf}_{\gamma }(f)\mid \dot{K}_{q}^{\alpha ,p}f_{\beta }^{s}\big\|%
=c\Big\|\Big(\sum\limits_{m\in \mathbb{Z}^{n}}2^{js}\inf_{j-\gamma ,m,\gamma
}(f)\chi _{j-\gamma ,m}\Big)_{j\geq \gamma }\mid \dot{K}_{q}^{\alpha
,p}(\ell _{\beta })\Big\|
\end{equation*}%
Define a sequence $\{\lambda _{i,k}\}_{i\in \mathbb{N}_{0},k\in \mathbb{Z}%
^{n}}$ by setting $\lambda _{i,k}=\inf_{y\in Q_{i,k}}|\widetilde{\varphi
_{i-\gamma }}\ast f(y)|$ and $\lambda _{0,k}=\inf_{y\in Q_{\gamma ,k}}|%
\widetilde{\Phi }\ast f(y)|$. We have%
\begin{equation*}
\inf_{j-\gamma ,m,\gamma }(f)=\sup_{h\in \mathbb{Z}^{n}}\{\lambda
_{j,h}:Q_{j,h}\cap Q_{j-\gamma ,m}\neq \emptyset \}
\end{equation*}%
and 
\begin{equation*}
\inf_{0,m,\gamma }(f)=\sup_{h\in \mathbb{Z}^{n}}\{\lambda _{0,h}:Q_{\gamma
,h}\cap Q_{0,m}\neq \emptyset \}.
\end{equation*}%
Let $h\in \mathbb{Z}^{n}$ with $Q_{j,h}\cap Q_{j-\gamma ,m}\neq \emptyset $
and $j\geq \gamma $. Then 
\begin{equation}
\lambda _{j,h}\leq c2^{\gamma d/r}\lambda _{j,z,r,\gamma }^{\ast }\text{, }%
j>\gamma \text{ \ \ and \ \ }\lambda _{0,h}\leq c2^{\gamma d/r}\lambda
_{0,z,r,\gamma }^{\ast }\text{, \ \ }j=\gamma   \label{est1lamda}
\end{equation}%
for any $z\in \mathbb{Z}^{n}$ with $Q_{j,z}\cap Q_{j-\gamma ,m}\neq
\emptyset $, where the constant $c>0$ does not depend on $j$, $h$ and $z$.
Indeed, we observe%
\begin{equation*}
\lambda _{j,h}=\frac{\lambda _{j,h}}{(1+2^{j}|2^{-j}h-2^{-j}z|)^{d}}%
(1+2^{j}|2^{-j}h-2^{-j}z|)^{d}.
\end{equation*}%
Let $x\in Q_{j,h}\cap Q_{j-\gamma ,m}$ and $y\in Q_{j,z}\cap Q_{j-\gamma ,m}$%
. We have%
\begin{equation*}
|2^{-j}h-2^{-j}z|\leq |2^{-j}h-x|+|x-y|+|y-2^{-j}z|\lesssim
2^{2-j}+2^{\gamma -j}.
\end{equation*}%
This implies \eqref{est1lamda}. Hence%
\begin{equation*}
\sum\limits_{m\in \mathbb{Z}^{n}}\inf_{j-\gamma ,m,\gamma }(f)\chi
_{j-\gamma ,m}\leq c\sum\limits_{k\in \mathbb{Z}^{n}}\lambda
_{j,k,r,d}^{\ast }\chi _{j,k},\text{ \ \ }j>\gamma 
\end{equation*}%
and%
\begin{equation*}
\sum\limits_{m\in \mathbb{Z}^{n}}\inf_{0,m,\gamma }(f)\chi _{0,m}\leq
c\sum\limits_{k\in \mathbb{Z}^{n}}\lambda _{0,k,r,d}^{\ast }\chi _{\gamma
,k},\text{ \ \ }j=\gamma ,
\end{equation*}%
with $r=\min (q,\frac{n}{\frac{n}{q}+\alpha },\beta )$ and $d>n$. Therefore,%
\begin{equation*}
\big\|\text{inf}_{\gamma }(f)\mid \dot{K}_{q}^{\alpha ,p}f_{\beta }^{s}\big\|%
\leq c\Big\|\Big(\sum\limits_{k\in \mathbb{Z}^{n}}2^{js}\lambda
_{j,k,r,d}^{\ast }\chi _{j,k}\Big)_{j\geq \gamma }\mid \dot{K}_{q}^{\alpha
,p}(\ell _{\beta })\Big\|.
\end{equation*}%
Notice that if $j=\gamma $ we replace $\lambda _{j,k,r,d}^{\ast }\chi _{j,k}$
by $\lambda _{0,k,r,d}^{\ast }\chi _{\gamma ,k}$. Applying Lemma \ref%
{lamda-equi} to estimate this term by%
\begin{equation*}
c\Big\|\Big(\sum\limits_{k\in \mathbb{Z}^{n}}2^{js}\lambda _{j,k}\chi _{j,k}%
\Big)_{j>\gamma }\mid \dot{K}_{q}^{\alpha ,p}(\ell _{\beta })\Big\|+c\Big\|%
\sum\limits_{k\in \mathbb{Z}^{n}}\lambda _{0,k}\chi _{\gamma ,k}\mid \dot{K}%
_{q}^{\alpha ,p}\Big\|,
\end{equation*}%
wich is bounded by%
\begin{equation*}
c\Big\|\Big(2^{js}\widetilde{\varphi _{j-\gamma }}\ast f\Big)_{j>\gamma
}\mid \dot{K}_{q}^{\alpha ,p}(\ell _{\beta })\Big\|+c\Big\|\widetilde{\Phi }%
\ast f\mid \dot{K}_{q}^{\alpha ,p}\Big\|.
\end{equation*}%
By Theorem \ref{loc-mean-char}, we obtain%
\begin{equation*}
\big\|\text{inf}_{\gamma }(f)\mid \dot{K}_{q}^{\alpha ,p}f_{\beta }^{s}\big\|%
\leq c\big\|f\mid \dot{K}_{q}^{\alpha ,p}F_{\beta }^{s}\big\|^{\prime }\leq c%
\big\|f\mid \dot{K}_{q}^{\alpha ,p}F_{\beta }^{s}\big\|,
\end{equation*}%
where we use the caracterization of Herz-type Triebel-Lizorkin spaces by
local means. Applying Lemma A.4 of \cite{FJ90}, see also Lemma 8.3 of \cite%
{Bow07}, we obtain 
\begin{equation*}
\text{inf}_{\gamma }(f)_{\min (q,\frac{n}{\frac{n}{q}+\alpha },\beta
),\gamma }^{\ast }\approx \sup (f)_{\min (q,\frac{n}{\frac{n}{q}+\alpha }%
,\beta ),\gamma }^{\ast }.
\end{equation*}%
Hence for $\gamma >0$ sufficiently large we obtain by applying Lemma \ref%
{lamda-equi}, 
\begin{equation*}
\big\|\text{inf}_{\gamma }(f)_{\min (q,\frac{n}{\frac{n}{q}+\alpha },\beta
),\gamma }^{\ast }\mid \dot{K}_{q}^{\alpha ,p}f_{\beta }^{s}\big\|\approx %
\big\|\text{inf}_{\gamma }(f)\mid \dot{K}_{q}^{\alpha ,p}f_{\beta }^{s}\big\|
\end{equation*}%
and 
\begin{equation*}
\big\|\sup (f)_{\min (q,\frac{n}{\frac{n}{q}+\alpha },\beta ),\gamma }^{\ast
}\mid \dot{K}_{q}^{\alpha ,p}f_{\beta }^{s}\big\|\approx \big\|\sup (f)\mid 
\dot{K}_{q}^{\alpha ,p}f_{\beta }^{s}\big\|
\end{equation*}%
for any $s\in \mathbb{R},0<p,q<\infty ,0<\beta \leq \infty $ and $\alpha
>-n/q$. Therefore, 
\begin{equation*}
\big\|\text{inf}_{\gamma }(f)\mid \dot{K}_{q}^{\alpha ,p}f_{\beta }^{s}\big\|%
\approx \big\|f\mid \dot{K}_{q}^{\alpha ,p}F_{\beta }^{s}\big\|\approx \big\|%
\sup (f)\mid \dot{K}_{q}^{\alpha ,p}f_{\beta }^{s}\big\|.
\end{equation*}%
Use these estimates and repeating the proof of Theorem 2.2 in \cite{FJ90} or
Theorem 2.1 in \cite{SiYY10} then complete the proof of Theorem \ref%
{phi-tran}. \ \ \rule{3mm}{3mm}

\begin{remark}
$\mathrm{From}$ $\mathrm{these}$ $\mathrm{to}$ $\mathrm{prove}$ $\mathrm{the}
$ $\mathrm{embeddings}$%
\begin{equation*}
\dot{K}_{q}^{\alpha _{2},r}F_{\theta }^{s_{2}}\hookrightarrow \dot{K}%
_{s}^{\alpha _{1},p}F_{\beta }^{s_{1}}
\end{equation*}%
$\mathrm{we}$ $\mathrm{need}$ $\mathrm{only}$ $\mathrm{to}$ $\mathrm{prove}$%
\begin{equation*}
\dot{K}_{q}^{\alpha _{2},r}f_{\theta }^{s_{2}}\hookrightarrow \dot{K}%
_{s}^{\alpha _{1},p}f_{\beta }^{s_{1}},
\end{equation*}%
$\mathrm{under}$ $\mathrm{the}$ $\mathrm{same}$ $\mathrm{restrictions}$ $%
\mathrm{on}$ $\mathrm{parameters}$ $s_{1},s_{2},\alpha _{1},\alpha
_{2},s,p,q,\beta ,r,\theta $.
\end{remark}

We end this section with one more lemma, which is basically a consequence of
Hardy's inequality in the sequence Lebesgue space $\ell _{q}$.

\begin{lemma}
\label{lq-inequality}\textit{Let }$0<a<1$\textit{\ and }$0<q\leq \infty $%
\textit{. Let }$\left\{ \varepsilon _{k}\right\} $\textit{\ be a sequences
of positive real numbers and denote }$\delta
_{k}=\sum_{j=0}^{k}a^{k-j}\varepsilon _{j}$ and $\eta
_{k}=\sum_{j=k}^{\infty }a^{j-k}\varepsilon _{j},k\in \mathbb{N}_{0}$. Then
there exists constant $c>0\ $\textit{depending only on }$a$\textit{\ and }$q$
such that%
\begin{equation*}
\Big(\sum\limits_{k=0}^{\infty }\delta _{k}^{q}\Big)^{1/q}+\Big(%
\sum\limits_{k=0}^{\infty }\eta _{k}^{q}\Big)^{1/q}\leq c\text{ }\Big(%
\sum\limits_{k=0}^{\infty }\varepsilon _{k}^{q}\Big)^{1/q}.
\end{equation*}
\end{lemma}

\section{Sobolev embeddings for $\dot{K}_{q}^{\protect\alpha ,p}F_{\protect%
\beta }^{s}$ spaces}

It is well-known that%
\begin{equation*}
F_{q,\infty }^{s_{2}}\hookrightarrow F_{s,\beta }^{s_{1}}
\end{equation*}%
if $s_{1}-n/s=s_{2}-n/q$, where $0<q<s<\infty $ and $0<\beta \leq \infty $
(see e.g. \cite[Theorem 2.7.1]{Triebel83}). In this section we generalize
these embeddings to Herz-type Triebel-Lizorkin spaces. We need the Sobolev
embeddings properties of the above sequence spaces.

\begin{theorem}
\label{embeddings3}\textit{Let }$\alpha _{1},\alpha _{2},s_{1},s_{2}\in 
\mathbb{R},0<s,q<\infty ,0<r\leq p<\infty ,0<\beta \leq \infty ,\alpha
_{1}>-n/s\ $\textit{and }$\alpha _{2}>-n/q$. \textit{We suppose that }%
\begin{equation}
s_{1}-n/s-\alpha _{1}=s_{2}-n/q-\alpha _{2}.  \label{newexp1}
\end{equation}

\noindent \textit{Let }$0<q<s<\infty $ and $\alpha _{2}\geq \alpha _{1}$ or $%
0<s\leq q<\infty $ and 
\begin{equation}
\alpha _{2}+n/q\geq \alpha _{1}+n/s.  \label{newexp2}
\end{equation}%
Then%
\begin{equation}
\dot{K}_{q}^{\alpha _{2},r}f_{\theta }^{s_{2}}\hookrightarrow \dot{K}%
_{s}^{\alpha _{1},p}f_{\beta }^{s_{1}},  \label{Sobolev-emb1}
\end{equation}%
where%
\begin{equation*}
\theta =\left\{ 
\begin{array}{ccc}
\beta & \text{if} & 0<s\leq q<\infty \text{ and }\alpha _{2}+n/q=\alpha
_{1}+n/s \\ 
\infty &  & \text{otherwise.}%
\end{array}%
\right.
\end{equation*}
\end{theorem}

\proof 
We would \ like to mention that this embedding was proved in \cite[Theorem
5.9]{Drihem13} under the restriction $\mathrm{max(}0,\frac{\alpha _{1}s}{q}%
\mathrm{)}\leq \alpha _{2}\leq \left( \alpha _{1}+\frac{n}{s}\right) \frac{r%
}{p}-\frac{n}{q}$. Here we use a different method to omit this condition.

\noindent \textit{Step 1.} Let us prove that $0<r\leq p<\infty $ is
necessary. In the calculations below we consider the 1-dimensional case for
simplicity. For any $v\in \mathbb{N}_{0}$ and $N\geq 2$, we put 
\begin{equation*}
\lambda _{v,m}^{N}=\left\{ 
\begin{array}{ccc}
2^{-(s_{1}-\frac{1}{s}-\alpha _{1})v}\sum_{i=2}^{N-2}\chi _{i}(2^{v-1}) & 
\text{if} & m=1 \\ 
0 &  & \text{otherwise,}%
\end{array}%
\right. 
\end{equation*}%
$\lambda ^{N}=\{\lambda _{v,m}^{N}:v\in \mathbb{N}_{0},m\in \mathbb{Z}\}$.
We have 
\begin{equation*}
\big\|\lambda ^{N}\mid \dot{K}_{s}^{\alpha _{1},p}f_{\beta }^{s_{1}}\big\|%
^{p}=\sum_{k=-\infty }^{\infty }2^{\alpha _{1}kp}\Big\|\Big(%
\sum_{v=0}^{\infty }\sum\limits_{m\in \mathbb{Z}}2^{vs_{1}\beta }|\lambda
_{v,m}^{N}|^{\beta }\chi _{v,m}\Big)^{1/\beta }\chi _{k}\Big\|_{s}^{p}.
\end{equation*}%
We can rewrite the last statement as follows: 
\begin{eqnarray*}
&&\sum_{k=1-N}^{0}2^{\alpha _{1}kp}\Big\|\Big(\sum_{v=2}^{N-2}2^{(\frac{1}{s}%
+\alpha _{1})v\beta }\chi _{v,1}\Big)^{1/\beta }\chi _{k}\Big\|_{s}^{p} \\
&=&\sum_{k=1-N}^{0}2^{\alpha _{1}kp}\big\|2^{(\frac{1}{s}+\alpha
_{1})(1-k)}\chi _{1-k,1}\big\|_{s}^{p}=c\text{ }N,
\end{eqnarray*}%
where the constant $c>0$ does not depend on $N$. Now%
\begin{equation*}
\big\|\lambda ^{N}\mid \dot{K}_{q}^{\alpha _{2},r}f_{\theta }^{s_{2}}\big\|%
^{r}=\sum_{k=-\infty }^{\infty }2^{\alpha _{2}kr}\Big\|\Big(%
\sum_{v=0}^{\infty }2^{vs_{2}\theta }|\lambda _{v,1}^{N}|^{\theta }\chi
_{v,1}\Big)^{1/\theta }\chi _{k}\Big\|_{q}^{r}.
\end{equation*}%
Again we can rewrite the last statement as follows:%
\begin{eqnarray*}
&&\sum_{k=1-N}^{0}2^{\alpha _{2}kr}\Big\|\Big(%
\sum_{v=2}^{N-2}2^{(s_{2}-s_{1}+\frac{1}{s}+\alpha _{1})v\theta }\chi _{v,1}%
\Big)^{1/\theta }\chi _{k}\Big\|_{q}^{r} \\
&=&\sum_{k=1-N}^{0}2^{\alpha _{2}kr}\Big\|2^{(s_{2}-s_{1}+\frac{1}{s}+\alpha
_{1})(1-k)}\chi _{1-k,1}\Big\|_{q}^{r}=c\text{ }N,
\end{eqnarray*}%
where the constant $c>0$ does not depend on $N$. If the embeddings %
\eqref{Sobolev-emb1} holds then for any $N\in \mathbb{N}$, $N^{\frac{1}{p}-%
\frac{1}{r}}\leq C$. Thus, we conclude that $0<r\leq p<\infty $ must
necessarily hold by letting $N\rightarrow +\infty $.

\textit{Step 2.} We consider the sufficiency of the conditions. First we
consider $0<q<s<\infty $ and $\alpha _{2}\geq \alpha _{1}$. In view of the
embedding $\ell _{r}\hookrightarrow \ell _{p}$, it is sufficient to prove
that 
\begin{equation*}
\dot{K}_{q}^{\alpha _{2},r}f_{\infty }^{s_{2}}\hookrightarrow \dot{K}%
_{s}^{\alpha _{1},r}f_{\beta }^{s_{1}}.
\end{equation*}%
By similarity, we only consider the case $\beta =1$. Let $\lambda \in \dot{K}%
_{q}^{\alpha _{2},r}f_{\infty }^{s_{2}}$. We have%
\begin{eqnarray}
\left\Vert \lambda \mid \dot{K}_{s}^{\alpha _{1},r}f_{1}^{s_{1}}\right\Vert 
&\leq &\Big(\sum\limits_{k=-\infty }^{c_{n}+1}2^{k\alpha _{1}r}\Big\|%
\sum_{v=0}^{\infty }2^{vs_{1}}\sum_{m\in \mathbb{Z}^{n}}\lambda _{v,m}\chi
_{v,m}\chi _{k}\Big\|_{s}^{r}\Big)^{1/r}  \notag \\
&&+\Big(\sum\limits_{k=c_{n}+2}^{\infty }2^{k\alpha _{1}r}\Big\|%
\sum_{v=0}^{\infty }2^{vs_{1}}\sum_{m\in \mathbb{Z}^{n}}\lambda _{v,m}\chi
_{v,m}\chi _{k}\Big\|_{s}^{r}\Big)^{1/r}.  \label{est 1}
\end{eqnarray}%
Here $c_{n}=1+[\log _{2}(2\sqrt{n}+1)]$. The first term can be estimated by%
\begin{eqnarray*}
&&c\Big(\sum\limits_{k=-\infty }^{c_{n}+1}2^{k\alpha _{1}r}\Big\|%
\sum_{v=0}^{c_{n}-k+1}2^{vs_{1}}\sum_{m\in \mathbb{Z}^{n}}\lambda _{v,m}\chi
_{v,m}\chi _{k}\Big\|_{s}^{r}\Big)^{1/r} \\
&&+c\Big(\sum\limits_{k=-\infty }^{c_{n}+1}2^{k\alpha _{1}r}\Big\|%
\sum_{v=c_{n}-k+2}^{\infty }2^{vs_{1}}\sum_{m\in \mathbb{Z}^{n}}\lambda
_{v,m}\chi _{v,m}\chi _{k}\Big\|_{s}^{r}\Big)^{1/r} \\
&=&I+II.
\end{eqnarray*}%
\textit{Estimation of }$I$\textit{.} Let $x\in C_{k}\cap Q_{v,m}$ and $y\in
Q_{v,m}$. We have $|x-y|\leq 2\sqrt{n}2^{-v}<2^{c_{n}-v}$ and from this it
follows that $|y|<2^{c_{n}-v}+2^{k}\leq 2^{c_{n}-v+2}$, which implies that $y
$ is located in some ball $B(0,2^{c_{n}-v+2})$. Then%
\begin{equation*}
|\lambda _{v,m}|^{t}=2^{nv}\int_{\mathbb{R}^{n}}|\lambda _{v,m}|^{t}\chi
_{v,m}(y)dy\leq 2^{nv}\int_{B\left( 0,2^{c_{n}-v+2}\right) }|\lambda
_{v,m}|^{t}\chi _{v,m}(y)dy,
\end{equation*}%
if $x\in C_{k}\cap Q_{v,m}$ and $t>0$. Therefore for any $x\in C_{k}$%
\begin{eqnarray*}
\sum_{m\in \mathbb{Z}^{n}}|\lambda _{v,m}|^{t}\chi _{v,m}(x) &\leq
&2^{nv}\int_{B\left( 0,2^{c_{n}-v+2}\right) }\sum_{m\in \mathbb{Z}%
^{n}}|\lambda _{v,m}|^{t}\chi _{v,m}(y)dy \\
&=&2^{nv}\Big\|\sum_{m\in \mathbb{Z}^{n}}|\lambda _{v,m}|\chi _{v,m}\chi
_{B\left( 0,2^{c_{n}-v+2}\right) }\Big\|_{t}^{t}.
\end{eqnarray*}%
Hence%
\begin{eqnarray*}
&&2^{\alpha _{1}k}\Big\|\sum_{v=0}^{c_{n}-k+1}2^{vs_{1}}\sum_{m\in \mathbb{Z}%
^{n}}\lambda _{v,m}\chi _{v,m}\chi _{k}\Big\|_{s} \\
&\leq &c\text{ }2^{(\alpha _{1}+\frac{n}{s})k}%
\sum_{v=0}^{c_{n}-k+1}2^{v(s_{1}+\frac{n}{t})}\Big\|\sum_{m\in \mathbb{Z}%
^{n}}|\lambda _{v,m}|\chi _{v,m}\chi _{B\left( 0,2^{c_{n}-v+2}\right) }\Big\|%
_{t}.
\end{eqnarray*}%
We may choose $t>0$ such that $\frac{1}{t}>\max (\frac{1}{q},\frac{1}{q}+%
\frac{\alpha _{2}}{n})$. Using \eqref{newexp1} and Lemma \ref{lq-inequality}
to estimate $I^{r}$ by%
\begin{eqnarray*}
&&c\sum_{v=0}^{\infty }2^{v(s_{2}-\frac{n}{q}-\alpha _{2}+\frac{n}{t})r}%
\Big\|\sum_{m\in \mathbb{Z}^{n}}|\lambda _{v,m}|\chi _{v,m}\chi _{B\left(
0,2^{c_{n}-v+2}\right) }\Big\|_{t}^{r} \\
&\leq &c\sum_{v=0}^{\infty }2^{v(s_{2}-\frac{n}{q}-\alpha _{2}+\frac{n}{t})r}%
\Big(\sum_{i\leq -v}\Big\|\sum_{m\in \mathbb{Z}^{n}}|\lambda _{v,m}|\chi
_{v,m}\chi _{i+c_{n}+2}\Big\|_{t}^{\sigma }\Big)^{r/\sigma } \\
&\leq &c\sum_{v=0}^{\infty }2^{v\frac{nr}{d}}\Big(\sum_{i\leq -v}2^{i\frac{%
n\sigma }{d}+\alpha _{2}\sigma i}\Big\|\sup_{j}\sum_{m\in \mathbb{Z}%
^{n}}2^{s_{2}j}|\lambda _{j,m}|\chi _{j,m}\chi _{i+c_{n}+2}\Big\|%
_{q}^{\sigma }\Big)^{r/\sigma },
\end{eqnarray*}%
by H\"{o}lder's inequality, with $\sigma =\min (1,t)$ and $\frac{n}{d}=\frac{%
n}{t}-\frac{n}{q}-\alpha _{2}$. Again, we apply Lemma \ref{lq-inequality} to
obtain%
\begin{equation*}
I^{r}\leq c\sum_{i=0}^{\infty }2^{-\alpha _{2}ir}\Big\|\sup_{j}\sum_{m\in 
\mathbb{Z}^{n}}2^{s_{2}j}|\lambda _{j,m}|\chi _{j,m}\chi _{-i+c_{n}+2}\Big\|%
_{q}^{r}\leq c\big\|\lambda \mid \dot{K}_{q}^{\alpha _{2},r}f_{\infty
}^{s_{2}}\big\|^{r}.
\end{equation*}%
\textit{Estimation of }$II$\textit{.} We see that it suffices to show that
for any $k\leq c_{n}+1$%
\begin{eqnarray*}
2^{k\alpha _{1}}\Big\|\sum_{v=c_{n}-k+2}^{\infty }2^{vs_{1}}\sum_{m\in 
\mathbb{Z}^{n}}\lambda _{v,m}\chi _{v,m}\chi _{k}\Big\|_{s} &\leq &C^{s/q}%
\text{ }2^{k\alpha _{2}}\Big\|\sup_{v\geq c_{n}-k+2}\sum_{m\in \mathbb{Z}%
^{n}}2^{vs_{2}}\lambda _{v,m}\chi _{v,m}\chi _{\widetilde{C}_{k}}\Big\|_{q}
\\
&=&\delta ,
\end{eqnarray*}%
where $\widetilde{C}_{k}=\{x\in \mathbb{R}^{n}:2^{k-2}<\left\vert
x\right\vert <2^{k+2}\}$ and $C=2\max (\frac{2^{\frac{1}{q}}}{1-2^{\frac{n}{s%
}-\frac{n}{q}}},\frac{2^{\frac{n}{s}}}{2^{\frac{n}{s}}-1})$. This claim can
be reformulated as showing that%
\begin{equation*}
\int_{C_{k}}2^{k\alpha _{1}s}\delta ^{-s}\Big(\sum_{v=c_{n}-k+2}^{\infty
}2^{vs_{1}}\sum_{m\in \mathbb{Z}^{n}}|\lambda _{v,m}|\chi _{v,m}(x)\Big)%
^{s}dx\leq 1.
\end{equation*}%
The left-hand side can be rewritten us 
\begin{eqnarray*}
&&\int_{C_{k}}\delta ^{-s}\Big(\sum_{v=c_{n}-k+2}^{\infty }2^{v(\frac{n}{s}-%
\frac{n}{q})+(\alpha _{1}-\alpha _{2})(v+k)+vs_{2}+k\alpha _{2}}\sum_{m\in 
\mathbb{Z}^{n}}|\lambda _{v,m}|\chi _{v,m}(x)\Big)^{s}dx \\
&\leq &\int_{C_{k}}\delta ^{-s}\Big(\sum_{v=c_{n}-k+2}^{\infty }2^{v(\frac{n%
}{s}-\frac{n}{q})+vs_{2}+k\alpha _{2}}\sum_{m\in \mathbb{Z}^{n}}|\lambda
_{v,m}|\chi _{v,m}(x)\Big)^{s}dx=T_{k},
\end{eqnarray*}%
since $\alpha _{2}\geq \alpha _{1}$. Let us prove that $T_{k}\leq 1$ for any 
$k\leq c_{n}+1$. Our estimate use partially some decomposition techniques
already used in \cite{Vybiral09}.

\noindent \textit{Case 1.} $\sup_{v\geq c_{n}-k+2,m}2^{vs_{2}+k\alpha
_{2}}|\lambda _{v,m}|\chi _{v,m}(x)\leq \delta $. In this case we obtain%
\begin{equation*}
T_{k}\leq C^{s}\int_{C_{k}}\Big(\delta ^{-1}\sup_{v\geq c_{n}-k+2}\sum_{m\in 
\mathbb{Z}^{n}}2^{vs_{2}+k\alpha _{2}}|\lambda _{v,m}|\chi _{v,m}(x)\Big)%
^{q}dx\leq 1.
\end{equation*}

\noindent \textit{Case 2.} $\sup_{v\geq c_{n}-k+2,m}2^{vs_{2}+k\alpha
_{2}}|\lambda _{v,m}|\chi _{v,m}(x)>\delta $. We can distinguish two cases
as follows:

$\bullet $ $\delta ^{-1}\sup_{v\geq c_{n}-k+2,m}2^{vs_{2}+k\alpha
_{2}}|\lambda _{v,m}|\chi _{v,m}(x)=\infty $, then there is nothing to prove.

$\bullet $ $\delta <\sup_{v\geq c_{n}-k+2,m}2^{vs_{2}+k\alpha _{2}}|\lambda
_{v,m}|\chi _{v,m}(x)<\infty $. Let $N\in \mathbb{N}$ be such that 
\begin{equation*}
2^{\frac{nN}{q}}<\delta ^{-1}\sup_{v\geq c_{n}-k+2}\sum_{m\in \mathbb{Z}%
^{n}}2^{vs_{2}+k\alpha _{2}}|\lambda _{v,m}|\chi _{v,m}(x)<2^{\frac{(N+1)n}{q%
}}.
\end{equation*}%
\textit{Subcase 2.1.} $k\geq c_{n}-N+2$. We split the sum over $v\geq
c_{n}-k+2$ into two parts, 
\begin{equation*}
\sum_{v=c_{n}-k+2}^{\infty }\cdot \cdot \cdot =\sum_{v=c_{n}-k+2}^{N}\cdot
\cdot \cdot +\sum_{v=N+1}^{\infty }\cdot \cdot \cdot .
\end{equation*}%
Let $x\in C_{k}\cap Q_{v,m}$ and $y\in Q_{v,m}$. We have $|x-y|\leq 2\sqrt{n}%
2^{-v}<2^{c_{n}-v}$ and from this it follows that $%
2^{k-2}<|y|<2^{c_{n}-v}+2^{k}<2^{k+2}$, which implies that $y$ is located in 
$\widetilde{C}_{k}$. Then%
\begin{equation*}
|\lambda _{v,m}|^{q}=2^{nv}\int_{\mathbb{R}^{n}}|\lambda _{v,m}|^{q}\chi
_{v,m}(y)dy\leq 2^{nv}\int_{\widetilde{C}_{k}}|\lambda _{v,m}|^{q}\chi
_{v,m}(y)dy,
\end{equation*}%
if $x\in C_{k}\cap Q_{v,m}$. But this immediately implies that%
\begin{eqnarray*}
\sum_{m\in \mathbb{Z}^{n}}|\lambda _{v,m}|^{q}\chi _{v,m}(x) &\leq
&2^{nv}\int_{\widetilde{C}_{k}}\sum_{m\in \mathbb{Z}^{n}}|\lambda
_{v,m}|^{q}\chi _{v,m}(y)dy \\
&=&2^{nv}\big\|\sum_{m\in \mathbb{Z}^{n}}|\lambda _{v,m}|\chi _{v,m}\chi _{%
\widetilde{C}_{k}}\big\|_{q}^{q}\leq 2^{(\frac{n}{q}-s_{2})qv-\alpha
_{2}qk}\delta ^{q}.
\end{eqnarray*}%
Therefore for any $x\in C_{k}$%
\begin{equation*}
\delta ^{-1}\sum_{v=c_{n}-k+2}^{N}2^{v(\frac{n}{s}-\frac{n}{q}%
)+vs_{2}+k\alpha _{2}}\sum_{m\in \mathbb{Z}^{n}}|\lambda _{v,m}|\chi
_{v,m}(x)\leq \sum_{v=c_{n}-k}^{N}2^{\frac{n}{s}v}\leq C\text{ }2^{\frac{n}{s%
}N}
\end{equation*}%
and 
\begin{eqnarray*}
&&\delta ^{-1}\sum_{v=N+1}^{\infty }2^{v(\frac{n}{s}-\frac{n}{q}%
)+vs_{2}+k\alpha _{2}}\sum_{m\in \mathbb{Z}^{n}}|\lambda _{v,m}|\chi
_{v,m}(x) \\
&=&2^{(\frac{n}{s}-\frac{n}{q})N}\delta ^{-1}\sum_{v=N+1}^{\infty }2^{(v-N)(%
\frac{n}{s}-\frac{n}{q})+vs_{2}+k\alpha _{2}}\sum_{m\in \mathbb{Z}%
^{n}}|\lambda _{v,m}|\chi _{v,m}(x) \\
&\leq &c\text{ }2^{(\frac{n}{s}-\frac{n}{q})N}\delta ^{-1}\sup_{v\geq
c_{n}-k+2}\sum_{m\in \mathbb{Z}^{n}}2^{vs_{2}+k\alpha _{2}}|\lambda
_{v,m}|\chi _{v,m}(x) \\
&\leq &C\text{ }2^{\frac{n}{s}N}.
\end{eqnarray*}%
Hence%
\begin{equation*}
T_{k}\leq C^{s}\int_{C_{k}}2^{nN}dx\leq C^{s}\int_{C_{k}}\Big(\delta
^{-1}\sup_{v\geq c_{n}-k+2}\sum_{m\in \mathbb{Z}^{n}}2^{vs_{2}+k\alpha
_{2}}|\lambda _{v,m}|\chi _{v,m}(x)\Big)^{q}dx\leq 1.
\end{equation*}%
\textit{Subcase 2.2.} $k<c_{n}-N+2$. We use the same of arguments as in
Subcase 2.1, in view of the fact that $\sum_{v=c_{n}-k+2}^{\infty }\cdot
\cdot \cdot \leq \sum_{v=N+1}^{\infty }\cdot \cdot \cdot $.

\noindent \textit{Estimate of \eqref{est 1}}. The arguments here are quite
similar to those used in the estimation of $II$. This complete the proof of
the first case.

Now we consider the case $0<s\leq q<\infty $ and $\alpha _{2}+n/q>\alpha
_{1}+n/s$. We only need to estimate the part $T_{k}$. H\"{o}lder's
inequality implies that%
\begin{eqnarray*}
T_{k} &\leq &\Big\|\delta ^{-1}\sum_{v=c_{n}-k+2}^{\infty }2^{(\frac{n}{s}-%
\frac{n}{q}+\alpha _{1}-\alpha _{2})(v+k)}2^{vs_{2}+k\alpha _{2}}\sum_{m\in 
\mathbb{Z}^{n}}|\lambda _{v,m}|\chi _{v,m}\chi _{k}\Big\|_{q}^{s} \\
&\leq &\Big\|\delta ^{-1}\sup_{v\geq c_{n}-k+2}\sum_{m\in \mathbb{Z}%
^{n}}2^{vs_{2}+k\alpha _{2}}|\lambda _{v,m}|\chi _{v,m}\chi _{k}\Big\|%
_{q}^{s}=C^{-s/q},
\end{eqnarray*}%
where the last inequality follows by the fact that $\alpha _{2}+n/q>\alpha
_{1}+n/s$. The remaining case can be easily solved. The proof is complete. \
\ \rule{3mm}{3mm}

As a corollary of Theorems \ref{phi-tran} and \ref{embeddings3}, we have the
following Sobolev embedding for $\dot{K}_{q}^{\alpha ,p}F_{\beta }^{s}$
spaces.

\begin{theorem}
\label{embeddings3.1}\textit{Let }$\alpha _{1},\alpha _{2},s_{1},s_{2}\in 
\mathbb{R},0<s,q<\infty ,0<r\leq p<\infty ,0<\beta \leq \infty ,\alpha
_{1}>-n/s\ $\textit{and }$\alpha _{2}>-n/q$. \textit{We suppose that }%
\begin{equation*}
s_{1}-n/s-\alpha _{1}=s_{2}-n/q-\alpha _{2}.
\end{equation*}%
\noindent \textit{Let }$0<q<s<\infty $ and $\alpha _{2}\geq \alpha _{1}$ or $%
0<s\leq q<\infty $ and 
\begin{equation*}
\alpha _{2}+n/q\geq \alpha _{1}+n/s.
\end{equation*}%
Then%
\begin{equation*}
\dot{K}_{q}^{\alpha _{2},r}F_{\theta }^{s_{2}}\hookrightarrow \dot{K}%
_{s}^{\alpha _{1},p}F_{\beta }^{s_{1}},
\end{equation*}%
where%
\begin{equation*}
\theta =\left\{ 
\begin{array}{ccc}
\beta & \text{if} & 0<s\leq q<\infty \text{ and }\alpha _{2}+n/q=\alpha
_{1}+n/s \\ 
\infty &  & \text{otherwise.}%
\end{array}%
\right.
\end{equation*}
\end{theorem}

\begin{remark}
$\mathrm{We}$ $\mathrm{would}$ $\mathrm{like}$ $\mathrm{to}$ $\mathrm{mention%
}$ $\mathrm{that}$ \eqref{newexp2} $\mathrm{and}$ $s_{1}-n/s-\alpha _{1}\leq
s_{2}-n/q-\alpha _{2}$ $\mathrm{are}$ $\mathrm{necessary,}$ $\mathrm{see}$ $%
\mathrm{\cite{Drihem13}}$.
\end{remark}

From Theorem\ \ref{embeddings3.1}\ and the fact that $\dot{K}%
_{s}^{0,s}F_{\beta }^{s_{1}}=F_{s,\beta }^{s_{1}}$\ we immediately arrive at
the following corollaries.

\begin{corollary}
\label{embeddings4}Let $s_{1},s_{2}\in \mathbb{R},0<s,q<\infty
,s_{1}-n/s=s_{2}-n/q-\alpha _{2}$, $0<r\leq s<\infty $ and $0<\beta \leq
\infty $. \textit{Let }$0<q<s<\infty $\ and $\alpha _{2}\geq 0$ or $0<s\leq
q<\infty $ and $\alpha _{2}+n/q\geq n/s$. Then%
\begin{equation*}
\dot{K}_{q}^{\alpha _{2},r}F_{\theta }^{s_{2}}\hookrightarrow F_{s,\beta
}^{s_{1}},
\end{equation*}%
where%
\begin{equation*}
\theta =\left\{ 
\begin{array}{ccc}
\beta & \text{if} & 0<s\leq q<\infty \text{ and }\alpha _{2}+n/q=n/s \\ 
\infty &  & \text{otherwise.}%
\end{array}%
\right.
\end{equation*}
\end{corollary}

\begin{corollary}
\label{embeddings5}Let $s_{1},s_{2}\in \mathbb{R},0<s,q<\infty
,s_{1}-n/s-\alpha _{1}=s_{2}-n/q,0<q\leq p<\infty $ and $0<\beta \leq \infty 
$. \textit{Let }$0<q<s<\infty $ and $\alpha _{1}\leq 0$ or $0<s\leq q<\infty 
$ and $n/q\geq \alpha _{1}+n/s$. Then%
\begin{equation*}
F_{q,\theta }^{s_{2}}\hookrightarrow \dot{K}_{s}^{\alpha _{1},p}F_{\beta
}^{s_{1}}.
\end{equation*}%
where%
\begin{equation*}
\theta =\left\{ 
\begin{array}{ccc}
\beta & \text{if} & 0<s\leq q<\infty \text{ and }n/q=\alpha _{1}+n/s \\ 
\infty &  & \text{otherwise.}%
\end{array}%
\right.
\end{equation*}
\end{corollary}

From the above corollaries and the fact that $\dot{K}_{q}^{\alpha
,r}F_{2}^{0}=\dot{K}_{q}^{\alpha ,r}$ for $1<r,q<\infty $ and $-\frac{n}{q}%
<\alpha <n-\frac{n}{q}$, see \cite{XuYang03} we obtain the following
embeddings between Herz and Triebel-Lizorkin spaces%
\begin{equation*}
\dot{K}_{q}^{\alpha _{2},r}\hookrightarrow F_{s,\beta }^{s_{1}},
\end{equation*}%
if $n/s-s_{1}=n/q+\alpha _{2}$, $1<r\leq s<\infty $, $0<\beta \leq \infty $,
and 
\begin{equation*}
1<q<s<\infty \ \ \ \text{and \ \ }0\leq \alpha _{2}<n-\frac{n}{q}
\end{equation*}%
or 
\begin{equation*}
1<s\leq q<\infty \text{ \ \ and \ \ }\frac{n}{s}-\frac{n}{q}<\alpha _{2}<n-%
\frac{n}{q}
\end{equation*}%
or 
\begin{equation*}
1<s\leq q<\infty ,\alpha _{2}=\frac{n}{s}-\frac{n}{q}\text{ \ \ and \ \ }%
\beta =2.
\end{equation*}%
Again we obtain%
\begin{equation*}
F_{q,\theta }^{s_{2}}\hookrightarrow \dot{K}_{s}^{\alpha _{1},p}
\end{equation*}%
holds if $n/s+\alpha _{1}=n/q-s_{2},0<\max (q,1)<p<\infty $ (or $%
1<q=p<\infty $), $0<\theta \leq \infty $ and 
\begin{equation*}
0<\max (q,1)<s<\infty \text{ \ \ and \ \ }-\frac{n}{s}<\alpha _{1}\leq 0.
\end{equation*}%
or%
\begin{equation*}
1<s\leq q<\infty \text{ \ \ and \ \ }-\frac{n}{s}<\alpha _{1}<n/q-n/s
\end{equation*}%
or%
\begin{equation*}
1<s\leq q<\infty \text{, }\alpha _{1}=n/q-n/s\text{ \ \ and \ \ }\theta =2.
\end{equation*}

From the Jawerth-Franke embeddings we have%
\begin{equation*}
F_{t,\infty }^{s_{3}}\hookrightarrow B_{q,t}^{s_{2}}\hookrightarrow
F_{s,\beta }^{s_{1}},
\end{equation*}%
if $s_{1},s_{2},s_{3}\in \mathbb{R},s_{1}-n/s=s_{2}-n/q=s_{3}-n/t,0<t<q<s<%
\infty $\ and $0<\beta \leq \infty $, see \cite[p. 60]{Triebel06}. Using our
results, we have the following useful consequences.

\begin{corollary}
Let $s_{1},s_{2},s_{3}\in \mathbb{R},0<s,q,t<\infty
,s_{1}-n/s=s_{2}-n/q=s_{3}-n/t$\ and $0<\beta \leq \infty $. Then%
\begin{equation*}
F_{t,\infty }^{s_{3}}\hookrightarrow \dot{K}_{q}^{0,s}F_{\infty
}^{s_{2}}\hookrightarrow F_{s,\beta }^{s_{1}},\quad 0<t\leq q<s<\infty .
\end{equation*}
\end{corollary}

To prove this it is sufficient to take in Corollary \ref{embeddings4}, $r=s$
and $\alpha _{2}=0$.\ However the desired embeddings are an immediate
consequence of the fact that 
\begin{equation*}
F_{t,\infty }^{s_{3}}\hookrightarrow F_{q,\infty }^{s_{2}}=\dot{K}%
_{q}^{0,q}F_{\infty }^{s_{2}}\hookrightarrow \dot{K}_{q}^{0,s}F_{\infty
}^{s_{2}}.
\end{equation*}

\begin{corollary}
Let $s_{1},s_{2}\in \mathbb{R},0<s,q<\infty ,s_{1}-n/s=s_{2}-n/q$\ and $%
0<\beta \leq \infty $. Then%
\begin{equation*}
F_{q,\infty }^{s_{2}}\hookrightarrow \dot{K}_{s}^{0,q}F_{\beta
}^{s_{1}}\hookrightarrow F_{s,\beta }^{s_{1}},\quad 0<q<s<\infty .
\end{equation*}
\end{corollary}

To prove this it is sufficient to take in Corollary \ref{embeddings5}, $p=q$
and $\alpha _{1}=0$. Then the desired embeddings are an immediate
consequence of the fact that 
\begin{equation*}
F_{q,\infty }^{s_{2}}\hookrightarrow \dot{K}_{s}^{0,q}F_{\beta
}^{s_{1}}\hookrightarrow \dot{K}_{s}^{0,s}F_{\beta }^{s_{1}}=F_{s,\beta
}^{s_{1}}.
\end{equation*}

\section{Applications}

In this section, we give a simple application of Theorems \ref{embeddings3}
and \ref{embeddings3.1}.

\begin{theorem}
\label{unco-basis}\textit{Let }$s\in \mathbb{R},0<p,q,\beta <\infty $\textit{%
\ and }$\alpha >-n/q$\textit{. }Then there exists a linear isomorphism $T$
which maps $\dot{K}_{q}^{\alpha ,p}F_{\beta }^{s}$ onto $\dot{K}_{q}^{\alpha
,p}f_{\beta }^{s}$. Moreover, there is an unconditional basis in $\dot{K}%
_{q}^{\alpha ,p}F_{\beta }^{s}$.
\end{theorem}

The mapping $T$ is generated by an appropriate wavelet system. A proof of
this theorem can be found in Xu \cite{Xu09} for the non-homogeneous
Herz-type Triebel-Lizorkin spaces and $\alpha >0$. This result is also true
for the spaces $\dot{K}_{q}^{\alpha ,p}F_{\beta }^{s}$, with $\alpha >-n/q$.
Indeed, the problem can be reduced to proof the $\dot{K}_{q}^{\alpha
,p}f_{\beta }^{s}$-version of Lemma 3.5 in Xu \cite{Xu09}. Therefore we need
to recall the definition of molecules.

\begin{definition}
\label{Atom-Def}Let $K,L\in \mathbb{N}_{0}$ and let $M>0$. A $K$-times
continuously differentiable function $a\in C^{K}(\mathbb{R}^{n})$ is called $%
[K,L,M]$-molecule concentrated in $Q_{v,m}$, if for some $v\in \mathbb{N}%
_{0} $ and $m\in \mathbb{Z}^{n}$

\begin{equation}
|D^{\alpha }a(x)|\leq 2^{v|\alpha |}(1+2^{v}|x-2^{-v}m|)^{-M}\text{,\quad
for\quad }0\leq |\alpha |\leq K,x\in \mathbb{R}^{n}  \label{diff-cond}
\end{equation}%
and if%
\begin{equation}
\int\limits_{\mathbb{R}^{n}}x^{\alpha }a(x)dx=0,\text{\quad for\quad }0\leq
|\alpha |<L\text{ and }v\geq 1.  \label{mom-cond}
\end{equation}
\end{definition}

If the molecule $a\ $is concentrated in $Q_{v,m}$, that means if it fulfills 
$\mathrm{\eqref{diff-cond}}$ and $\mathrm{\eqref{mom-cond}}$, then we will
denote it by $a_{vm}$. For $v=0$ or $L=0$ there are no moment\ conditions\ $%
\mathrm{\eqref{mom-cond}}$ required.

Now, we prove the $\dot{K}_{q}^{\alpha ,p}f_{\beta }^{s}$-version of Lemma
3.5 in Xu \cite{Xu09}.

\begin{lemma}
\textit{Let }$s\in \mathbb{R},0<p,q<\infty ,0<\beta \leq \infty $\textit{\
and }$\alpha >-n/q$\textit{. }Furthermore, let $K,L\in \mathbb{N}_{0}$ and
let $M>0$ with%
\begin{equation*}
L>n(\frac{1}{\min (1,q,\beta )}-1)-1-s,\text{\quad }K\text{ arbitrary and }M%
\text{ large enough}.
\end{equation*}%
If $a_{vm}$ are $\left[ K,L,M\right] $-molecules concentrated in $Q_{v,m}$
and $\lambda =\{\lambda _{v,m}\in \mathbb{C}:v\in \mathbb{N}_{0},m\in 
\mathbb{Z}^{n}\}\in \dot{K}_{q}^{\alpha ,p}f_{\beta }^{s}$, then the sum%
\begin{equation}
\sum_{v=0}^{\infty }\sum_{m\in \mathbb{Z}^{n}}\lambda _{v,m}a_{vm}
\label{new-rep}
\end{equation}%
converges in $\mathcal{S}^{\prime }(\mathbb{R}^{n})$.
\end{lemma}

\textbf{Proof.} We use the arguments of \cite{Triebel06}, see also \cite{K10}%
. Let $\varphi \in \mathcal{S(}\mathbb{R}^{n})$. We get from the moment
conditions $\mathrm{\eqref{mom-cond}}$ for fixed $v\in \mathbb{N}_{0}$%
\begin{eqnarray*}
&&\int\limits_{\mathbb{R}^{n}}\sum\limits_{m\in \mathbb{Z}^{n}}\lambda
_{v,m}a_{vm}(y)\varphi (y)dy \\
&=&\int\limits_{\mathbb{R}^{n}}\sum\limits_{m\in \mathbb{Z}^{n}}\lambda
_{v,m}2^{-v(L+1)}a_{vm}(y)\Big(\varphi (y)-\sum\limits_{\left\vert \beta
\right\vert <L}(y-2^{-v}m)^{\beta }\frac{D^{\beta }\varphi (2^{-v}m)}{\beta !%
}\Big)2^{v(L+1)}dy \\
&=&\sum\limits_{i=-\infty }^{\infty }\int\limits_{C_{i}}\cdot \cdot \cdot dy,
\end{eqnarray*}%
where $C_{i}=\{y\in \mathbb{R}^{n}:2^{i-1}\leq \left\vert y\right\vert
<2^{i}\}$ for any $i\in \mathbb{Z}$. Let us estimate the sum $%
\sum\limits_{i=-\infty }^{0}\cdot \cdot \cdot $. We use the Taylor expansion
of $\varphi $ up to order $L-1$\ with respect to the off-points $2^{-v}m$,
we obtain%
\begin{equation*}
\varphi (y)-\sum\limits_{\left\vert \beta \right\vert <L}(y-2^{-v}m)^{\beta }%
\frac{D^{\beta }\varphi (2^{-v}m)}{\beta !}=\sum\limits_{\left\vert \beta
\right\vert =L}(y-2^{-v}m)^{\beta }\frac{D^{\beta }\varphi (\xi )}{\beta !},
\end{equation*}%
with $\xi $ on the line segment joining $y$ and $2^{-v}m$. Since $%
1+\left\vert y\right\vert \leq \left( 1+\left\vert \xi \right\vert \right)
\left( 1+\left\vert y-2^{-v}m\right\vert \right) $, we estimate%
\begin{eqnarray*}
\Big|\sum\limits_{\left\vert \beta \right\vert =L}(y-2^{-v}m)^{\beta }\frac{%
D^{\beta }\varphi (\xi )}{\beta !}\Big| &\leq &\left( 1+\left\vert
y-2^{-v}m\right\vert \right) ^{L}\sum\limits_{\left\vert \beta \right\vert
=L}\frac{|D^{\beta }\varphi (\xi )|}{\beta !} \\
&\leq &\left( 1+\left\vert y-2^{-v}m\right\vert \right) ^{L}(1+\left\vert
\xi \right\vert )^{-S}\left\Vert \varphi \right\Vert _{S,L} \\
&\leq &c\left( 1+\left\vert y\right\vert \right) ^{-S}\left( 1+\left\vert
y-2^{-v}m\right\vert \right) ^{L+S},
\end{eqnarray*}%
where $S>0$ is at our disposal. Let $0<t<\min (1,q)=1+q-\frac{q}{\min (1,q)}$
and $h=s+\frac{n}{q}(t-1)$ be such that $n(1-\frac{1}{\min (1,q)})+s>h>-1-L$%
. Since $a_{vm}$ are $\left[ K,L,M\right] $-molecules, then $%
2^{-v(L+1)}\left\vert a_{vm}(y)\right\vert \leq 2^{hv}2^{-v(L+1+h)}\left(
1+2^{v}\left\vert y-2^{-v}m\right\vert \right) ^{-M}$. Therefore, The sum $%
\sum\limits_{i=-\infty }^{0}\cdot \cdot \cdot $ can be estimated by%
\begin{equation}
c\text{ }2^{-v(L+1+h)}\sum\limits_{i=-\infty
}^{0}\int\limits_{C_{i}}\sum\limits_{m\in \mathbb{Z}^{n}}2^{hv}\left\vert
\lambda _{v,m}\right\vert \left( 1+2^{v}\left\vert y-2^{-v}m\right\vert
\right) ^{L+S-M}(1+\left\vert y\right\vert )^{-S}dy.  \label{Max-estimate}
\end{equation}%
Since $M$ can be taken large enough, by Lemma 4 in \cite{K10} we obtain%
\begin{equation*}
\sum\limits_{m\in \mathbb{Z}^{n}}\left\vert \lambda _{v,m}\right\vert \left(
1+2^{v}\left\vert y-2^{-v}m\right\vert \right) ^{L+S-M}\leq c\mathcal{M}\Big(%
\sum\limits_{m\in \mathbb{Z}^{n}}\left\vert \lambda _{v,m}\right\vert \chi
_{v,m}\Big)(y)
\end{equation*}%
for any $y\in C_{i}\cap Q_{v,l}$ with $l\in \mathbb{Z}^{n}$. We split $S$
into $R+T$ with $R+\alpha <0$ and $T$ large enough such that $T>\max (-R,%
\frac{n(q-t)}{q})$. Then $\mathrm{\eqref{Max-estimate}}$ is bounded by%
\begin{equation*}
c\ 2^{-v(L+1+h)}\sum\limits_{i=-\infty }^{0}2^{-iR}\int\limits_{C_{i}}%
\mathcal{M}\Big(\sum\limits_{m\in \mathbb{Z}^{n}}2^{vh}\left\vert \lambda
_{v,m}\right\vert \chi _{v,m}\Big)(y)(1+\left\vert y\right\vert )^{-T}dy.
\end{equation*}%
Since we have in addition the factor $(1+\left\vert y\right\vert )^{-T}$, it
follows by H\"{o}lder's inequality that this expression is bounded by%
\begin{eqnarray*}
&&c\ 2^{-v(L+1+h)}\sum\limits_{i=-\infty }^{0}2^{-iR}\Big\|\mathcal{M}\Big(%
\sum\limits_{m\in \mathbb{Z}^{n}}2^{hv}\left\vert \lambda _{v,m}\right\vert
\chi _{v,m}\Big)\chi _{i}\Big\|_{q/t} \\
&\leq &c\text{ }2^{-v(L+1+h)}\sum\limits_{i=-\infty }^{0}2^{-i(\alpha +R)}%
\Big\|\sum\limits_{m\in \mathbb{Z}^{n}}2^{hv}\left\vert \lambda
_{v,m}\right\vert \chi _{v,m}|\dot{K}_{q/t}^{\alpha ,\infty }\Big\| \\
&\leq &c\text{ }2^{-v(L+1+h)}\Big\|\lambda |\dot{K}_{q/t}^{\alpha
,p}f_{\infty }^{h}\Big\|,
\end{eqnarray*}%
where the first inequality follows by the boundedness of the
Hardy-Littlewood maximal operator $\mathcal{M}$ on $\dot{K}_{q/t}^{\alpha
,\infty }$. Using a combination of the arguments used above, the sum $%
\sum\limits_{i=1}^{\infty }\cdot \cdot \cdot $ can be estimated by $c$ $%
2^{-v(L+1+h)}\big\|\lambda |\dot{K}_{q/t}^{\alpha ,p}f_{\infty }^{h}\big\|$.
Since $L+1+h>0$, the convergence of $\mathrm{\eqref{new-rep}}$ is now clear
by the embeddings 
\begin{equation*}
\dot{K}_{q}^{\alpha ,p}f_{\infty }^{s}\hookrightarrow \dot{K}_{q/t}^{\alpha
,p}f_{\infty }^{h},
\end{equation*}%
see Theorem \ref{embeddings3}. The proof is completed. \ \ \rule{3mm}{3mm}

Let $w$ denote a positive, locally integrable function and $0<p<\infty $.
Then the weighted Lebesgue space $L^{p}(\mathbb{R}^{n},w)$ contains all
measurable functions such that 
\begin{equation*}
\big\|f\mid L^{p}(\mathbb{R}^{n},w)\big\|=\Big(\int_{\mathbb{R}%
^{n}}\left\vert f(x)\right\vert ^{p}w(x)dx\Big)^{1/p}<\infty .
\end{equation*}%
For $\varrho \in \lbrack 1,\infty )$ we denote by $\mathcal{A}_{\varrho }$
the Muckenhoupt class of weights, and $\mathcal{A}_{\infty }=\cup _{\varrho
\geq 1}\mathcal{A}_{\varrho }$. We refer to \cite{GR85} for the general
properties of these classes. Let $w\in \mathcal{A}_{\infty }$, $s\in \mathbb{%
R}$, $0<\beta \leq \infty $ and $0<p<\infty $. We define weighted
Triebel-Lizorkin spaces $F_{p,q}^{s}(\mathbb{R}^{n},w)$ to be the set of all
distributions $f\in \mathcal{S}^{\prime }(\mathbb{R}^{n})$ such that%
\begin{equation*}
\big\|f\mid F_{p,\beta }^{s}(\mathbb{R}^{n},w)\big\|=\Big\|\Big(%
\sum\limits_{j=0}^{\infty }2^{js\beta }\left\vert \mathcal{F}^{-1}\varphi
_{j}\ast f\right\vert ^{\beta }\Big)^{1/\beta }\mid L^{p}(\mathbb{R}^{n},w)%
\Big\|
\end{equation*}%
is finite. In the limiting case $q=\infty $ the usual modification is
required. The spaces $F_{p,\beta }^{s}(\mathbb{R}^{n},w)=F_{p,\beta }^{s}(w)$
are independent of the particular choice of the smooth dyadic resolution of
unity $\{\varphi _{j}\}_{j\in \mathbb{N}_{0}}$ appearing in their
definitions. They are quasi-Banach spaces (Banach spaces for $p,q\geq 1$),
and 
\begin{equation*}
\mathcal{S}(\mathbb{R}^{n})\hookrightarrow F_{p,\beta
}^{s}(w)\hookrightarrow \mathcal{S}^{\prime }(\mathbb{R}^{n}).
\end{equation*}%
Moreover, for $w\equiv 1\in \mathcal{A}_{\infty }$ we obtain the usual
(unweighted) Triebel-Lizorkin spaces. Let $w_{\gamma }$ be a power weight,
i.e., $w_{\gamma }(x)=|x|^{\gamma }$ with $\gamma >-n$. Then in view of the
fact that $L^{p}=\dot{K}_{p}^{0,p}$, we have%
\begin{equation*}
\big\|f\mid F_{p,\beta }^{s}(w_{\gamma })\big\|\approx \big\|f\mid \dot{K}%
_{p}^{\frac{\gamma }{p},p}F_{\beta }^{s}\big\|\text{.}
\end{equation*}

Applying Corollary \ref{embeddings3.1} in some particular cases yields the
following embeddings, see for the case of .

\begin{corollary}
\textit{Let }$s_{1},s_{2}\in \mathbb{R}$, $0<q<s<\infty $, $0<\beta \leq
\infty $ and $w_{\gamma _{1}}(x)=|x|^{\gamma _{1}}$, $w_{\gamma
_{2}}(x)=|x|^{\gamma _{2}}$, with $\gamma _{1}>-n\ $\textit{and }$\gamma
_{2}>-n$. \textit{We suppose that } 
\begin{equation*}
s_{1}-\frac{n+\gamma _{1}}{s}=s_{2}-\frac{n+\gamma _{2}}{q}
\end{equation*}%
and%
\begin{equation*}
\gamma _{2}/q\geq \gamma _{1}/s.
\end{equation*}
Then%
\begin{equation*}
F_{q,\infty }^{s_{2}}(w_{\gamma _{2}})\hookrightarrow F_{s,\beta
}^{s_{1}}(w_{\gamma _{1}}).
\end{equation*}
\end{corollary}

\begin{remark}
$\mathrm{We}$ $\mathrm{refer}$ $\mathrm{the}$ $\mathrm{reader}$ $\mathrm{to}$
$\mathrm{the}$ $\mathrm{recent}$ $\mathrm{paper}$ $\mathrm{\cite{HarLkr08}}$ 
$\mathrm{for}$ $\mathrm{further}$ $\mathrm{results}$ $\mathrm{about}$ $%
\mathrm{Sobolev}$ $\mathrm{embeddings}$ $\mathrm{for}$ $\mathrm{weighted}$ $%
\mathrm{spaces}$ $\mathrm{of}$ $\mathrm{Besov}$ $\mathrm{type}$ $\mathrm{%
where}$ $\mathrm{the}$ $\mathrm{weight}$ $\mathrm{belongs}$ $\mathrm{to}$ $%
\mathrm{some}$ $\mathrm{Muckenhoupt}$ $\mathcal{A}_{\varrho }$ $\mathrm{%
class.}$ $\mathrm{Notice}$ $\mathrm{that}$ $\mathrm{this}$ $\mathrm{results}$
$\mathrm{are}$ $\mathrm{given}$ $\mathrm{in\ \cite[Theorem \ 1.2]{MV12}}$ $%
\mathrm{but}$ $\mathrm{under}$ $\mathrm{the}$ $\mathrm{restrictions}$ $%
1<q<s<\infty \mathrm{,}$ $1\leq \beta \leq \infty \mathrm{.}$
\end{remark}





\begin{thebibliography}{99}
\bibitem{Bow07} M. Bownik, Anisotropic Triebel-Lizorkin Spaces with Doubling
Measures. \textit{The Journal of Geometric Analysis}. \textbf{17}(3) (2007),
337-424.

\bibitem{Drihem13} {\normalsize D. Drihem, Embeddings properties on
Herz-type Besov and Triebel-Lizorkin spaces. \emph{Math. Ineq and Appl.} 
\textbf{16}(2) (2013), 439-460.}

\bibitem{FJ90} {\normalsize M. Frazier, B. Jawerth, A discrete transform and
decomposition of distribution spaces. \emph{J. Funct. Anal.} \textbf{93}(1)
(1990), 34-170.}

\bibitem{GR85} {\normalsize J. Garc\.{\i}a-Cuerva,J.L. Rubio de Francia, 
\emph{Weighted norm inequalities and related topics}. In: North-Holland
Mathematics Studies, Vol. 116, North-Holland, Amsterdam, 1985.}

\bibitem{K10} {\normalsize H. Kempka, Atomic, molecular and wavelet
decomposition of 2-microlocal Besov and Triebel-Lizorkin spaces with
variable integrability, \textit{Funct Approx}, \textbf{43} (1010), 171-208.}

\bibitem{HarLkr08} {\normalsize D. D. Haroske, L. Skrzypczak, Entropy and
approximation numbers of embeddings of function spaces with Muckenhoupt
weights. I. \emph{Rev. Mat. Complut.} \textbf{21}(1) (2008),135-177.}

\bibitem{HerYang99} {\normalsize E. Hernandez, D. Yang, Interpolation of
Herz-type Hardy spaces and applications. \emph{Math. Nachr.} \textbf{42}
(1998), 564-581.}

\bibitem{LiYang96} {\normalsize X. Li, D. Yang, Boundedness of some
sublinear operators on Herz spaces. \emph{Illinois. J. Math.} \textbf{40 }%
(1996), 484-501.}

\bibitem{LuYang1.95} {\normalsize S. Lu, D. Yang, The local versions of $%
H^{p}(\mathbb{R}^{n})$\ spaces at the origin. \emph{Studia. Math.} \textbf{%
116} (1995), 103-131.}

\bibitem{LuYang2.95} {\normalsize S. Lu, D. Yang, The decomposition of the
weighted Herz spaces on $\mathbb{R}^{n}$\ and its applications. \emph{Sci.
in. China (Ser.A).} \textbf{38} (1995), 147-158.}

\bibitem{LuYang97} {\normalsize S. Lu, D. Yang, Herz-type Sobolev and Bessel
potential spaces and their applications. \emph{Sci in China (Ser.A).} 
\textbf{40} (1997), 113-129.}

\bibitem{MV12} {\normalsize M. Meyries, M.C. Veraar, Sharp embedding results
for spaces of smooth functions with power weights. \emph{Studia. Math.} 
\textbf{208}(3) (2012), 257-293.}

\bibitem{TD00} {\normalsize L. Tang, D. Yang, Boundedness of vector-valued
operators on weighted Herz spaces. \emph{Approx. Th. \& its Appl.} \textbf{16%
} (2000), 58-70.}

\bibitem{Triebel83} {\normalsize H. Triebel, \emph{Theory of function spaces}%
. Basel: Birkh\"{a}user, 1983.}

\bibitem{Triebel92} {\normalsize H. Triebel, \emph{Theory of function spaces
II}. Basel: Birkh\"{a}user, 1992.}

\bibitem{Triebel06} {\normalsize H. Triebel, Fractals and spectra. Birkh\"{a}%
user, Basel 1997.}

\bibitem{T11} Y. Tsutsui, The Navier-Stokes equations and weak Herz spaces, 
\textit{Advances in Differential Equations}. \textbf{16} (2011) 1049-1085.

\bibitem{Vybiral09} {\normalsize J. Vyb\'{\i}ral, Sobolev and Jawerth
embeddings for spaces with variable smoothness and integrability.\emph{\
Ann. Acad. Sci. Fenn. Math.} \textbf{34}(2) (2009), 529-544.}

\bibitem{XuYang03} {\normalsize J. Xu, D. Yang, Applications of Herz-type
Triebel-Lizorkin spaces. \emph{Acta. Math. Sci (Ser. B).} \textbf{23}
(2003), 328-338.}

\bibitem{XuYang05} {\normalsize J. Xu, D. Yang, Herz-type Triebel-Lizorkin
spaces, I.\emph{\ Acta. Math. Sci (English Ed.)}. \textbf{21}(3) (2005),
643-654.}

\bibitem{Xu03} {\normalsize J. Xu, Some properties on Herz-type Besov spaces
(in chinese). \emph{J. Hunan. Univ (Natural Sci).} \textbf{30}(5) (2003),
75-78.}

\bibitem{Xu04} {\normalsize J. Xu, Pointwise multipliers of Herz-type Besov
spaces and their applications.\emph{\ Appl. Math. }\textbf{17}(1) (2004),
115-121.}

\bibitem{Xu05} {\normalsize J. Xu. Equivalent norms of Herz-type Besov and
Triebel-Lizorkin spaces. \emph{J. Funct. Spaces. Appl.} \textbf{3} (2005),
17-31.}

\bibitem{Xu09} {\normalsize J. Xu, Decompositions of non-homogeneous
Herz-type Besov and Triebel-Lizorkin spaces. \emph{Sci. China. Math. }}%
\textbf{57}(2) (2014), 315-331.

\bibitem{SiYY10} {\normalsize W. Yuan, W. Sickel, D. Yang, Morrey and
Campanato Meet Besov, Lizorkin and Triebel. Lecture Notes in Mathematics,
vol. 2005, Springer-Verlag, Berlin 2010.}
\end{thebibliography}
\end{document}